\documentclass[reqno, 12pt]{amsart}
\usepackage{amsfonts}
\usepackage{amsmath,amsthm,amscd, amssymb}
\usepackage{amsthm}

\title[Tian-Yau-Zelditch expansion]{Generalized asymptotic expansions of Tian-Yau-Zelditch (Announcements)}
\author{Chiung-ju Liu and Zhiqin Lu}
\date{March 15, 2007}

\thanks{
The first author is supported by the NSC grant NSC 982115M002007
in Taiwan. The second author is partially supported by the NSF
grant DMS 0904653.} \email[Zhiqin Lu, Department of Mathematics,
University of California, Irvine, CA 92697]{zlu@uci.edu}
\email[Chiung-ju Liu, Taida Institute for Mathematical Sciences,
Taiwan]{cjliu4@ntu.edu.tw}
\newtheorem{theorem}{Theorem}[section]
\newtheorem*{thm}{Theorem}

\newtheorem{cor}{Corollary}[section]

\theoremstyle{remark}
\newtheorem{rem}{Remark}[section]

\font\strange=msbm10

\newfont{\fnt}{cmr10 scaled 550}

\newcommand{\R}{{{\mathchoice  {\hbox{$\textstyle{\text{\strange R}}$}}
{\hbox{$\textstyle{\text{\strange R}}$}} {\hbox{$\scriptstyle
N\kern-0.3em  R$}} {\hbox{$\scriptscriptstyle  R\kern-0.2em
R$}}}}}

\newcommand{\Z}{{{\mathchoice  {\hbox{$\textstyle{\text{\strange Z}}$}}
{\hbox{$\textstyle{\text{\strange Z}}$}} {\hbox{$\scriptstyle
Z\kern-0.3em  Z$}} {\hbox{$\scriptscriptstyle  Z\kern-0.2em
Z$}}}}}

\newcommand{\N}{{{\mathchoice  {\hbox{$\textstyle{\text{\strange N}}$}}
{\hbox{$\textstyle{\text{\strange N}}$}} {\hbox{$\scriptstyle
N\kern-0.3em  N$}} {\hbox{$\scriptscriptstyle  N\kern-0.2em
N$}}}}}

\usepackage{color}

\renewcommand{\phi}{\varphi}

\usepackage{amsmath,amsthm,amscd}

\begin{document}
\maketitle
\section{Introduction}
Let $M$ be an $n$-dimensional projective algebraic manifold in
certain projective space $\mathbb{CP}^{N}$. The hyperplane line
bundle of $\mathbb{CP}^{N}$ restricts to an ample line bundle $L$
on $M$, which is called a polarization of $M$. A K\"ahler metric
$g$ is called a polarized metric,  if the corresponding K\"ahler
form represents the first Chern class $c_{1}(L)$ of $L$ in
$H^{2}(M,\mathbb{Z})$. Given any polarized K\"ahler metric $g$,
there is a Hermitian metric $h$ on $L$ whose Ricci form is equal
to $\omega_{g}$. For each positive integer $m>0$, the Hermitian
metric $h_L$ induces the Hermitian metric $h_{L}^{m}$ on $L^{m}$.
Let $(E,h_{E})$ be a Hermitian vector bundle of rank $r$ with a
Hermitian metric $h_{E}$. Consider the space
$\Gamma(M,L^{m}\otimes E)$ of all holomorphic sections for large
$m$. For $U,V\in\Gamma(M,L^{m}\otimes E)$, the pointwise and the
$L^{2}$ inner products are defined  as
\[\langle U(x),V(x)\rangle_{h_{L}^{m}\otimes h_{E}}\] and
\numberwithin{equation}{section}
\begin{equation}(U,V)=\int_{M}\langle
U(x),V(x)\rangle_{h_{L}^{m}\otimes
h_{E}}dV_{g},\label{1.1}\end{equation}
respectively,
where
$dV_{g}=\frac{\omega_{g}^{n}}{n!}$ is the volume form of $g$. Let
$\{S_{1},\cdots,S_{d}\}$ be an orthonormal basis of
$\Gamma(M,L^{m}\otimes E)$ with respect to \eqref{1.1}, where
$d=d(m)=\dim_{\mathbb{C}}\Gamma(M,L^{m}\otimes E)$. For any $x\in
M$, define a matrix $S=S(x)$ by
\[
S=\bigg(\langle S_{i}, S_{j}\rangle_{h_{L}^{m}\otimes h_{E}}\bigg).
\]
For any positive integer $b$, define
\begin{equation}\sigma_{b}\equiv {\rm tr}(S^{b}).\end{equation}
The value $\sigma_{b}$ is independent of the choice of the
orthonormal basis
because under different basis of $\Gamma(M,L^m\otimes E)$, the matrices $S$ are similar.
Moreover, $S$ is
diagonalizable. Since $E$ is of rank $r$, there exists a unitary
matrix $Q$ and a diagonal matrix $D$ such that $Q^{\ast}SQ=D$,
where
\begin{equation}
D_{ij}=\left\{%
\begin{array}{ll}
    \lambda_{i}\delta_{ij}, & \hbox{if $1\leq i,j \leq r$;} \\
    0, & \hbox{otherwise.} \\
\end{array}%
\right.    \label{1.2}
      \end{equation}

Furthermore, there exists an orthonormal basis $\{T_{i}\}_{i=1}^{d}$
such that
\begin{equation}\label{1.2.1}
\sigma_{b}={\rm tr}(S^{b})={\rm
tr}(D^{b})=\sum_{i=1}^{d}\|T_{i}(x)\|^{2b}_{h_{L}^{m}\otimes
h_{E}}.\end{equation}

From \eqref{1.2}, we have
\[\sigma_{b}=\sum_{i=1}^{r}\lambda_{i}^{b},\]
where $\lambda_{i}$ are the nonzero eigenvalues of $S$. For $b>r$,
$\sigma_{b}$ can be written as a polynomial of
$\sigma_{1},\cdots,\sigma_{r}$. Hence we only need to compute
$\sigma_{1},\cdots,\sigma_{r}$.

The asymptotic behavior of $\sigma_{b}$ plays a very important
rule in K\"ahler-Einstein geometry. In the case of $b=1$, Zelditch
\cite{zel} and Catlin \cite{cat} independently proved the existence
of an asymptotic expansion (Tian-Yau-Zeldtich expansion) of the
Szeg\"o kernel. In the break through paper of Donaldson \cite{skd1},
using the expansion, he was able to prove the stability for the
manifold admitting constant scalar curvature.

The result of Zeldtich and Catlin is stated as follows:
\begin{thm}[Zelditch, Catlin]\label{thm1.1}
Let $M$ be a compact complex manifold of dimension $n$ (over
$\mathbb{C}$) and let $(L,h)\rightarrow M$ be a positive Hermitian
holomorphic line bundle. Let $x$ be a point of $M$. Let $g$ be the
K\"ahler metric on $M$ corresponding to the K\"ahler form
$\omega_g=Ric(h)$. For each $m\in \N$, $h$ induces a Hermitian
metric $h_m$ on $L^m$. Let $\{S_1^m,\cdots,S_{d_m}^m\}$ be any
orthonormal basis of $H^0(M,L^m)$, $d_m=\dim H^0(M,L^m)$, with
respect to the inner product~\eqref{1.1}.
Then there is an asymptotic expansion:
\begin{equation}\label{fud1}
\sigma_1 = \sum_{i=1}^{d_m}||S_i^m(x)||_{h_m}^2 \sim
a_0(x)m^n+a_1(x)m^{n-1}+a_2(x)m^{n-2}+\cdots
\end{equation}
for certain smooth coefficients $a_j(x)$ with $a_0=1$. More
precisely, for any $k$
\[
||\sum_{i=1}^{d_m}||S_i^m(x)||_{h_m}^2
-\sum_{k=0}^Na_j(x)m^{n-k}||_{C^\mu}\leq C_{N,\mu}m^{n-N-1},
\]
where $C_{N,\mu}$ depends on $N,\mu$ and the manifold $M$.
\end{thm}
In \cite{lu10}, Lu proved that each coefficient $a_{j}(x)$ is a
polynomial of the curvature and its covariant derivatives. In
particular, $a_{1}(x)=\frac 12\rho(x)$ is half of the scalar
curvature of the K\"ahler manifold. All polynomials $a_{j}(x)$ can
be represented by a polynomial of the curvature and its
derivatives . Moreover, Lu and Tian \cite[Theorem 3.1]{lt-1}
proved that the leading term of $a_{j}$ is $C\triangle^{j-1}\rho$,
where $\rho$ is the scalar curvature and $C=C(j,n)$ is a constant.

In this paper, we establish an asymptotic expansion for
$\sigma_{b}$. Note that both Zelditch and Catlin used Szeg\"o krnel
or Bergman kernel in their proofs. Their methods, however, do not
apply to  the case $b>1$. To establish the expansion, we go back to
peak section estimates. Using the peak section method in
\cite{t5}, we are able to get the expansion of $\sigma_{b}$, which
generalizes the result of Zelditch and Catlin. Our result is:

\begin{theorem}\label{thm}
Let $M$ be a compact complex manifold of dimension $n$,
$(L,h_{L})\rightarrow M$ a positive Hermitian holomorphic line
bundle and $(E,h_{E})$ a Hermitian vector bundle of rank $r$. Let
$g$ be the K\"ahler metric on $M$ corresponding to the K\"ahler form
$\omega_{g}=Ric(h_{L})$. Let $\Gamma(M,L^{m}\otimes E)$ be the space
of all holomorphic global sections of $L^{m}\otimes E$, and let  $\{T_{1},\cdots,T_{d}\}$  be an orthonormal basis  of
$\Gamma(M,L^{m}\otimes E)$. Let
\begin{equation}\sigma_{b}=\sum_{i=1}^{d}\big\|T_{i}(x)\big\|_{h_{L}^{m}\otimes
h_{E}}^{2b}.\label{1.4}\end{equation} Then for $m$ big enough, there
exists an asymptotic expansion
\begin{equation}
\sigma_{b}(x)\sim
a_{0}(x)m^{bn}+a_{1}(x)m^{bn-1}+\cdots\label{1.5}
\end{equation}
for certain smooth coefficients $a_{j}(x)$. The expansion is in the
sense that
\begin{equation*}
\|\sigma_{b}-\sum_{k=0}^{N}a_{k}m^{bn-k}\|_{C^{\mu}}\leq
C(\mu,N,M)m^{bn-N-1}\end{equation*}
for positive integers $N,\mu$ and a constant $C(N,\mu, M)$ depending only on $N$, $\mu$ and the manifold.
\end{theorem}
The second main result of this paper focuses on compact complex
manifolds with analytic K\"ahler metrics. It is well-known that
the Tian-Yau-Zeldtich expansion does not converge in general. Even if it is
convergent, it may not converge to $\sigma_{1}$. We proved that in
the case when the metric is analytic, the optimal result may
achieve: The asymptotic expansion is convergent and the limit
approaches $\sigma_{b}$ faster than any other polynomials.
\begin{theorem}\label{thm1.2}
With the notations as in the above theorem, suppose that the
Hermitian metrics $h_{L}$  and $h_E$ are  real analytic at a fixed point $x$. Then
for $m$ big enough, the expansion
\begin{equation}
\sigma_{b}(x)\sim a_{0}(x)m^{bn}+a_{1}(x)m^{bn-1}+\cdots\label{1.6}
\end{equation}
is convergent for certain smooth coefficients $a_{j}(x)$. There is a $\delta>0$ such that the coefficient $a_{j}(x)$ satisfies
\[|a_{j}(x)|< \frac{C}{\delta^{j}}\] for some constant $C$.
Moreover, the expansion is convergent in the sense
\begin{equation*}
\|\sigma_{b}-\sum_{k=0}^{N}a_{k}m^{bn-k}\|_{C^{\mu}}\leq
Cm^{bn}(\delta m)^{-N-1}\end{equation*} for a constant $C(\mu,\delta)$
which only depends on $\mu$.
\end{theorem}
Theorem~\ref{thm1.2} gives that
\begin{cor}With the notations as in the above theorem, the limit of the series
\begin{equation*}
\lim_{N\rightarrow \infty}\sum_{k=0}^{N}a_{k}m^{bn-k}
\end{equation*}
exists.
\end{cor}
In fact, we prove a little bit more in Theorem \ref{thm1.2}.
\begin{cor}
With the notations as in the above theorem, we have
\begin{equation*}
\|\sigma_{b}-\sum_{k=0}^{N}a_{k}m^{bn-k}\|_{C^{\mu}}\leq
Ce^{-(\log m)^{2}}.\end{equation*}
\end{cor}
\begin{proof}
Choose $N=[\log m]$ to be the integer part of $\log m$. Then
\begin{equation*}
\|\sigma_{b}-\sum_{k=0}^{N}a_{k}m^{bn-k}\|_{C^{\mu}}\leq C
m^{bn-\log m-1}.\end{equation*} On the other hand,
\begin{equation*}
\|\sum_{k=N+1}^{\infty}a_{k}m^{bn-k}\|_{C^{\mu}}\leq
Cm^{bn-N-1}.\end{equation*} Thus \begin{equation*}
\|\sigma_{b}-\sum_{k=0}^{\infty}a_{k}m^{bn-k}\|_{C^{\mu}}\leq C
m^{bn-N}\leq Cm^{bn}e^{-(\log m)^{2}}.\end{equation*}
\end{proof}
 More
precisely, we have the following result
\begin{equation}
\|\sigma_{b}-\sum_{k=0}^{\infty}a_{k}m^{bn-k}\|_{C^{\mu}}\leq
Ce^{-\varepsilon(\log m)^{2}}.\label{1.8}\end{equation}

The above result \eqref{1.8} was only known in very special cases
before. In Liu \cite{cliu-1}, she proved the case for $b=1$ on a
smooth Riemann surface with constant curvature. On a planar domain
with Poincar\'e metric, Engli\u s proved the same result.

We have multiple definitions for $O(\frac{1}{m^{k+1}})$ through
this paper. In the case that $h_{L}$ and $h_E$ are $C^{\infty}$ as
the assumption in Theorem \ref{thm}, it denotes a quantity
dominated by $C/m^{k+1}$ with the constant $C$ depending only on
$k$ and the geometry of $M$. In the case that $h_{L}$ and $h_E$
are analytic as the assumption in Theorem \ref{thm1.2}, it denotes
a pure constant.

The last part of this paper, we compute the coefficient of the new
expansion.
\begin{theorem}\label{thm1.3}With the same notation as in Theorem \ref{thm}, each coefficient
$a_{j}(x)$ is a homogeneous polynomial of the curvature and its
derivatives at $x$. In particular,
\begin{align*}\begin{split}
a_{0}&=r\\
a_{1}&= \frac{1}{2}br\rho+\rho_E,
\end{split}
\end{align*}
and the leading term for $a_{k}$ for $k\geq 2$ is
\[\frac{brk}{(k+1)!}\triangle
^{k-1}\rho+\frac{1}{k!}\triangle^{k-1}\rho_E ,\] where $\rho$ is
the scalar curvature of $M$, $\rho_E$ is the scalar curvature of
$E$, and $\triangle$ is the Laplace operator of $M$.
\end{theorem}

\begin{rem} Apart from the generality of the above results,
the $C^0$-estimate of the Tian-Yau-Zelditch expansion, when $b\neq 1$,  was essentially known to~\cite{lu10}.
Thus, technically the proofs of the above theorems are on the $C^\mu$ estimates.
We developed a theorem on the variation of $K$-coordinates to achieve this goal.
This is a research announcement.
Details of the proofs will follow.
\end{rem}

\bibliographystyle{abbrv}

\bibliography{pub,unp}

\def\cprime{$'$} \def\cprime{$'$} \def\cprime{$'$}
\begin{thebibliography}{1}

\bibitem{cat}
D.~Catlin.
\newblock The {B}ergman kernel and a theorem of {T}ian.
\newblock In {\em Analysis and geometry in several complex variables (Katata,
  1997)}, Trends Math., pages 1--23. Birkh\"auser Boston, Boston, MA, 1999.

\bibitem{skd1}
S.~K. Donaldson.
\newblock Scalar curvature and projective embeddings. {I}.
\newblock {\em J. Differential Geom.}, 59(3):479--522, 2001.

\bibitem{cliu-1}
C.~J.~Liu.
\newblock {The asymptotic Tian-Yau-Zelditch expansion on Riemann surfaces with
  Constant Curvature}.
\newblock {\em (to be appeared in Taiwanese J. Math.)}, arXiv:0710.1347v3.

\bibitem{lu10}
Z.~Lu.
\newblock On the lower order terms of the asymptotic expansion of
  {T}ian-{Y}au-{Z}elditch.
\newblock {\em Amer. J. Math.}, 122(2):235--273, 2000.

\bibitem{lt-1}
Z.~Lu and G.~Tian.
\newblock The log term of the {S}zeg{\H o} kernel.
\newblock {\em Duke Math. J.}, 125(2):351--387, 2004.

\bibitem{t5}
G.~Tian.
\newblock On a set of polarized {K}\"ahler metrics on algebraic manifolds.
\newblock {\em J. Differential Geom.}, 32(1):99--130, 1990.

\bibitem{zel}
S.~Zelditch.
\newblock Szeg{\H o} kernels and a theorem of {T}ian.
\newblock {\em Internat. Math. Res. Notices}, no.(6):317--331, 1998.

\end{thebibliography}
\end{document}